\documentclass[reprint,a4paper,11pt,reqno]{amsart}

\usepackage{amsmath}    
\usepackage{amsfonts}
\usepackage{amscd}
\usepackage[latin2]{inputenc}
\usepackage{t1enc}
\usepackage[mathscr]{eucal}
\usepackage{indentfirst}
\usepackage{graphicx}
\usepackage{graphics}
\usepackage{pict2e}
\usepackage{fancyhdr}
\numberwithin{equation}{section}
\usepackage[margin=2.9cm]{geometry}
\usepackage{epstopdf} 
\usepackage{hyperref}
\usepackage{inputenc}

\theoremstyle{plain}
\newtheorem{Th}{Theorem}
\newtheorem{Lemma}[Th]{Lemma}
\newtheorem{Cor}[Th]{Corollary}

\theoremstyle{definition}

\newtheorem{?}[Th]{Problem}

\hypersetup{
    colorlinks=true,
    linkcolor=blue,
    citecolor=blue,
    filecolor=magenta,      
    urlcolor=cyan,
}

\begin{document}

\pagestyle{fancy}
\lhead{}
\chead{}
\rhead{\thepage}
\cfoot{}
\lfoot{}
\rfoot{}
\renewcommand{\headrule}{}

\title{On certain finite and infinite sums of inverse tangents} 

\author{Martin Nicholson} 

\begin{abstract}  
      An identity is proved connecting two finite sums of inverse tangents. This identity is discretized version of Jacobi's imaginary transformation for the modular angle from the theory of elliptic functions. Some other related identities are discussed.
\end{abstract}


\clearpage\maketitle
\thispagestyle{empty}
\vspace{-20pt}
\section{Introduction}

Sums of inverse tangents have attracted a lot of attention. For example, the following sums of inverse tangents can be calculated in closed form:
\begin{equation}\label{glaisher}
    \sum_{n=0}^\infty\arctan\frac{2}{(2n+1)^2}=\frac{\pi}{2},
\end{equation}
\begin{equation}\label{hr}
    \sum_{n=0}^\infty(-1)^{n+1}\arctan\frac{1}{F_{2n}}=\arctan\frac{\sqrt{5}-1}{2},
\end{equation}
\begin{equation}\label{bragg}
    \sum_{n=0}^\infty\arctan\frac{\sinh x}{\cosh nx}=\frac{3\pi}{4}-\arctan e^x.
\end{equation}
(\ref{glaisher}) is a classic sum evaluated first by Glaisher in \cite{glaisher}. The sum (\ref{hr}), where $F_n$ is $n$-th Fibonacci number, was calculated by Hoggatt and Ruegels \cite{hoggatt}. The sum (\ref{bragg}) was noted in \cite{bragg}. See \cite{berndt2} for further references and a brief summary of research in this direction.

All summations of the type (\ref{glaisher}) and (\ref{hr}) seem to be based on two methods: the telescopic principle, and the method of zeroes, as was noted in \cite{glasser}.

Even earlier, in his studies of elliptic functions, Jacobi proved identity of which he wrote in his treatise on elliptic functions ``one is obliged to rank among the most elegant formulas'' \cite{folkmar},\cite{ww}:
\begin{equation}\label{modular_angle}
    \frac{1}{4}\arcsin k=\arctan q^{1/2}-\arctan q^{3/2}+\arctan q^{5/2}-\ldots
\end{equation}
Here $q=e^{-\pi K'/K}$, $K$ is the complete elliptic integral of the first kind with modulus $k$
$$
K=K(k)=\int_0^{\pi/2}\frac{d\varphi}{\sqrt{1-k^2\sin^2\varphi}},
$$
$K'=K(k')$ with $k'=\sqrt{1-k^2}$ being the complementary modulus. The quantity $\arcsin k$ is called modular angle. Together with the obvious relation
$$
\arcsin k+\arcsin k'=\frac{\pi}{2},
$$
this implies
\begin{equation}\label{modular}
    \sum_{n=1}^{\infty}\chi(n)\arctan e^{-\alpha n}+\sum_{n=1}^{\infty}\chi(n)\arctan e^{-\beta n}=\frac{\pi}{8}, \qquad \alpha\beta=\frac{\pi^2}{4},
\end{equation}
where $\chi(n)=\sin\frac{\pi n}{2}$ is Dirichlet character modulo $4$ (\cite{berndt}, ch.14, entry 15). (\ref{modular}) is Jacobi's imaginary transformation for the modular angle.

Another arctan series related to elliptic functions was found in an unpublished manuscript by B. Cais \cite{cais}:
\begin{equation}\label{modular2}
    \sum_{n=1}^{\infty} \left(\frac{n}{3}\right)\arctan\frac{\sqrt{3}}{1+2 e^{\alpha n}}+\sum _{n=1}^{\infty} \left(\frac{n}{3}\right) \arctan\frac{\sqrt{3}}{1+2 e^{\beta n}}=\frac{\pi}{18}, \qquad \alpha\beta=\frac{4\pi^2}{9},
\end{equation}
where $\left(\frac{j}{3}\right) = \frac{2}{\sqrt{3}}\sin\frac{2\pi j}{3}$ is Legendre symbol modulo $3$.

The focus of this paper will be two reciprocal identities for finite sums of inverse tangents, Theorems \ref{th1} and \ref{th2} below and another transformation formula with two continuous parameters, Theorem \ref{th3}. We give two proofs of Theorem \ref{th1} in sections \ref{proof11}, \ref{proof12}. Theorems \ref{th2} and \ref{th3} are proved in sections \ref{proof2} and \ref{proof3}, respectively. In section \ref{discretedirichlet} we mention another transformation formula for sum of two finite reciprocal sums related to solution of Dirichlet problem on a finite rectangular grid.

\begin{Th}\label{th1} Let $n,m\in\mathbb{N}_0$ and $\alpha\beta=1$, $\alpha>0$. Then
\begin{align}\label{reciprocal_sum2}
\nonumber\sum_{|j|\le n}(-1)^{n+j}\arctan& {\biggl(\!\sqrt{1+\alpha^2\cos^2\!\tfrac{\pi  j}{2n+1}}-\alpha\cos\tfrac{\pi  j}{2n+1}\!\biggr)^{2 m+1}}\\&+\sum_{|k|\le m}(-1)^{m+k}\arctan\!{\biggl(\!\sqrt{1+\beta^2\cos ^2\!\tfrac{\pi  k}{2m+1}}-\beta\cos\tfrac{\pi k}{2m+1}\!\biggr)^{2 n+1}}=\frac{\pi}{4}.
\end{align}
\end{Th}
\noindent Note that when $n=m$ and $\alpha=1$ both sums in (\ref{reciprocal_sum2}) are equal and we get a closed form summation:
\begin{Cor} For $n\in\mathbb{N}_0$$\mathrm{:}$
\begin{equation}
    \sum_{|j|\le n}(-1)^{n+j}\arctan {\biggl(\!\sqrt{1+\cos^2\!\tfrac{\pi  j}{2n+1}}-\cos\tfrac{\pi  j}{2n+1}\!\biggr)^{2 n+1}}=\frac{\pi}{8}.
\end{equation}
\end{Cor}
\noindent It is instructive to write (\ref{reciprocal_sum2}) in another form by shifting the summation variable and simple rearrangement of terms
\begin{align*}\label{reciprocal_sum}
\nonumber\sum_{j=1}^{2n}\chi_4(j)&\arctan {\biggl(\!\sqrt{1+\alpha^2\sin^2\!\tfrac{\pi  j}{4n+2}}-\alpha\sin\tfrac{\pi  j}{4n+2}\!\biggr)^{2 m+1}}\\&+\sum_{k=1}^{2m}\chi_4(k)\arctan\!{\biggl(\!\sqrt{1+\beta^2\sin ^2\!\tfrac{\pi  j}{4m+2}}-\beta\sin\tfrac{\pi  j}{4m+2}\!\biggr)^{2 n+1}}\\&=\frac{\pi}{8}-\frac12(-1)^n\arctan\left(\sqrt{1+\alpha^2}-\alpha\right)^{2 m+1}-\frac12(-1)^m\arctan\left(\sqrt{1+\beta^2}-\beta\right)^{2 n+1}.
\end{align*}
From this form of (\ref{reciprocal_sum2}), it is evident that letting $n=m\to\infty$ and redefining $\alpha$ and $\beta$, one recovers (\ref{modular}). Thus, (\ref{reciprocal_sum2}) is discretized version of (\ref{modular}). Our proof is completely elementary and provides an elementary proof of the modular relation (\ref{modular}).

\begin{Th}\label{th2}
Let $n$ and $m$ be positive odd numbers and $\alpha\beta=1$. Then
\begin{equation}\label{reciprocalsum3}
\sum _{j=1}^{3n/2} \left(\frac{j}{3}\right) \arctan\frac{\sqrt{3}}{1+2\left(\cfrac{\alpha+\tan\frac{\pi j}{3n}}{\alpha-\tan\frac{\pi j}{3n}}\right)^{\!\! m}}+\sum _{k=1}^{3m/2} \left(\frac{k}{3}\right) \arctan\frac{\sqrt{3}}{1+2\left(\cfrac{\beta+\tan\frac{\pi k}{3m}}{\beta-\tan\frac{\pi k}{3m}}\right)^{\!\! n}}=-\frac{\pi}{6},
\end{equation}
where $\left(\frac{j}{3}\right)$ is Legendre symbol modulo $3$.
\end{Th}
\noindent If $n=m$ and $\alpha=1$ the two sums in Theorem \ref{th2} are equal and we get closed form summation
\begin{Cor} For an odd positive integer $n$$\mathrm{:}$
\begin{equation}\label{corollary2}
    \sum _{j=1}^{3n/2} \left(\frac{j}{3}\right) \arctan\frac{\sqrt{3}}{1+2\cot^n\biggl(\cfrac{\pi }{4}-\cfrac{\pi  j}{3 n}\biggr)}=-\frac{\pi}{12}.
\end{equation}
\end{Cor}
\noindent As an illustration of (\ref{corollary2}) note the case $n=3$:
$$
\arctan\frac{\sqrt{3}}{1+2\cot^3\frac{\pi}{36}}-
\arctan\frac{\sqrt{3}}{1+2\cot^3\frac{5\pi}{36}}-
\arctan\frac{\sqrt{3}}{1+2\cot^3\frac{29\pi}{36}}=\frac{\pi }{12}.
$$

Although (\ref{reciprocalsum3}) has a structure similar to (\ref{modular2}) it is not clear if (\ref{modular2}) can be derived from (\ref{reciprocalsum3}) as a limiting case.
However, by combining the limiting case of Theorem \ref{th2} with (\ref{modular2}) one can find the transformation formula for another infinite arctan series:
\begin{equation}\label{modular3}
    \sum_{n=0}^{\infty} \left(\frac{n-1}{3}\right)\arctan\frac{\sqrt{3}}{1-2 e^{\alpha (2n+1)}}+\sum _{n=0}^{\infty} \left(\frac{n-1}{3}\right) \arctan\frac{\sqrt{3}}{1-2 e^{\beta (2n+1)}}=\frac{2\pi}{9}, \qquad \alpha\beta=\frac{\pi^2}{9}.
\end{equation}
More generally one has as a consequence of imaginary transform for theta functions \cite{ww} (see Theorem \ref{th3} below for definition of $s(x)$)
$$
\sum _{j=-\infty}^\infty\!\! s(j)\arctan\frac{\sin 2 \theta}{e^{2 \alpha (\pi |j|+\varphi s(j) )}-\cos 2 \theta}+\sum _{k=-\infty}^\infty\!\! s(k)\arctan\frac{\sin 2 \varphi}{e^{2 \beta (\pi |k|+\theta s(k) )}-\cos 2 \varphi}=\frac{2 }{\pi }\left(\frac{\pi }{2}-\theta \right)\left(\frac{\pi }{2}-\varphi \right)\!,
$$
which suggests the following generalization of Theorem \ref{th2} with two additional continuous parameters.

\begin{Th}\label{th3} Let $n$ and $m$ be positive odd integers, $\alpha\beta=1$ $(\alpha>0)$, and $\theta,\varphi\in (0,\pi/2)$. Define function $s(x)=\begin{cases} \phantom{-}1,~ x\ge 0\\ -1,~x<0 \end{cases}
$, which differs from the $\mathrm{sgn}$ function only at $0$ where it takes value $1$ instead of value $0$. Then
\begin{align}\label{parameter}
    \nonumber \sum _{|j|\le \frac{n-1}{2}}s(j) \arctan&\frac{\sin 2 \theta }{\left(\cfrac{\alpha+\tan\frac{\varphi+\pi j}{n}}{\alpha-\tan\frac{\varphi+\pi  j}{n}}\right)^{m s(j)}-\cos2\theta}\\&+\sum _{|k|\le \frac{m-1}{2}}s(k) \arctan\frac{\sin 2 \varphi }{\left(\cfrac{\beta+\tan\frac{\theta+\pi k}{m}}{\beta-\tan\frac{\theta+\pi  k}{m}}\right)^{n s(k)}-\cos2\varphi}=\frac{\pi}{2}-\theta-\varphi.
\end{align}
\end{Th}
To see that Theorem \ref{th2} is a particular case of Theorem \ref{th3} let $\theta=\varphi=\pi/3$ in \ref{parameter}. The first sum in \ref{parameter} becomes
$$
\sum _{j=0}^{\frac{n-1}{2}}\arctan\frac{\sqrt{3} }{1+2\left(\cfrac{\alpha+\tan\frac{\pi(1+3j)}{3n}}{\alpha-\tan\frac{\pi(1+3j)}{3n}}\right)^{m}}-\sum_{j=1}^{\frac{n-1}{2}}\arctan\frac{\sqrt{3} }{1+2\left(\cfrac{\alpha-\tan\frac{\pi(1-3j)}{3n}}{\alpha+\tan\frac{\pi(1-3j)}{3n}}\right)^{m}}.
$$
After shifting the summation index of the second sum in this expression, one obtains
using the definition of Legendre symbol
\begin{align*}
\sum _{j=0}^{\frac{n-1}{2}}\arctan\frac{\sqrt{3} }{1+2\left(\cfrac{\alpha+\tan\frac{\pi(1+3j)}{3n}}{\alpha-\tan\frac{\pi(1+3j)}{3n}}\right)^{m}}&-\sum_{j=0}^{\frac{n-3}{2}}\arctan\frac{\sqrt{3} }{1+2\left(\cfrac{\alpha+\tan\frac{\pi(2+3j)}{3n}}{\alpha-\tan\frac{\pi(2+3j)}{3n}}\right)^{m}}\\&=\sum _{j=1}^{\frac{3n-1}{2}}\left(\frac{j}{3}\right)\arctan\frac{\sqrt{3} }{1+2\left(\cfrac{\alpha+\tan\frac{\pi j}{3n}}{\alpha-\tan\frac{\pi j}{3n}}\right)^{m}}.
\end{align*}

When $m=n$, $\varphi=\theta$, $\alpha=\beta=1$ in Theorem \ref{th3} one gets the summation formula
\begin{Cor} Let $n$ be a positive odd integer, $\theta\in (0,\pi/2)$, and the function $s(x)$ defined as in Theorem \ref{th3}. Then
\begin{equation}
     \sum _{|j|\le \frac{n-1}{2}}s(j) \arctan\frac{\sin 2 \theta }{\left\{\tan\left(\frac{\pi}{4}+\frac{\theta+\pi  j}{n}\right)\right\}^{n s(j)}-\cos2\theta}=\frac{\pi}{4}-\theta.
\end{equation}
\end{Cor}
As an illustration of the corollary above note the case $n=3$
\begin{align}
    \nonumber\arctan\frac{\sin 2 \theta }{\cos 2 \theta +\cot ^3\left(\frac{\pi }{12}+\frac{\theta }{3}\right)}&-\arctan\frac{\sin 2 \theta }{\cos 2 \theta +\cot ^3\left(\frac{\pi }{12}-\frac{\theta }{3}\right)}\\&+\arctan\frac{\sin 2 \theta }{\cos 2 \theta +\cot ^3\left(\frac{\theta }{3}-\frac{\pi }{4}\right)}=\theta-\frac{\pi}{4}.
\end{align}

Technically, one could generalize Theorem \ref{th1} too, but the resulting identity is not nice. We give it here for illustration purposes only without proof. Assuming $\theta,\varphi\in (0,\pi/2)$, $\alpha=1/\beta>0$, $n,m\in\mathbb{N}_0$ one has
\begin{align*}
{\mathrm{Re}}\Bigg[\sum_{|j|\le n}(-1)^{n+j}&\arctan\left\{ {\biggl(\!\sqrt{1+\alpha^2\cos^2\!\tfrac{\theta+\pi  j}{2n+1}}-\alpha\cos\tfrac{\theta+\pi  j}{2n+1}\!\biggr)^{2 m+1}}e^{i\varphi}\right\}\\&+\sum_{|k|\le m}(-1)^{m+k}\arctan\left\{\!{\biggl(\!\sqrt{1+\beta^2\cos ^2\!\tfrac{\varphi+\pi  k}{2m+1}}-\beta\cos\tfrac{\varphi+\pi k}{2m+1}\!\biggr)^{2 n+1}}e^{i\theta}\right\}\Bigg]=\frac{\pi}{4}.
\end{align*}

\section{First proof of Theorem \ref{th1}}\label{proof11}
We break the proof into a series of lemmas.
\begin{Lemma}\label{lemma1} The following identity holds for $\alpha>0$, $n,m\in\mathbb{N}_0$ and $j\in\mathbb{Z}$
$$
2\arctan {\biggl(\!\sqrt{1+\alpha^2\cos^2\!\tfrac{\pi  j}{2n+1}}-\alpha\cos\tfrac{\pi  j}{2n+1}\!\biggr)^{2 m+1}}=\frac{\pi}{2}-\arctan\left(\sinh(2m+1)\alpha_j\right)
$$
where $\alpha_j$ is the positive solution of $\sinh\alpha_j=\alpha \cos\frac{\pi j}{2n+1}$.
\end{Lemma}

\noindent{\it{Proof.}} By denoting $s=2m+1$ for brevity we obtain
\begin{align*}
    2\arctan {\biggl(\!\sqrt{1+\alpha^2\cos^2\!\tfrac{\pi  j}{2n+1}}-\alpha\cos\tfrac{\pi j}{2n+1}\!\biggr)^{2m+1}}&=2\arctan\left(\cosh\alpha_j-\sinh\alpha_j\right)^s\\
    &=2\arctan e^{-s\alpha_j}\\
    &=\frac{\pi}{2}-\left(\arctan e^{s\alpha_j}-\arctan e^{-s\alpha_j}\right)\\
    &=\frac{\pi}{2}-\arctan\frac{e^{s\alpha_j}-e^{-s\alpha_j}}{2}.
\end{align*}
Since ${\displaystyle{\frac{e^{x}-e^{-x}}{2}=\sinh x}}$ the proof is complete. \qed

\begin{Lemma} For $n,m\in\mathbb{N}_0$, $j\in\mathbb{Z}$, and $\alpha_j$ as was defined in the previous lemma, one has
$$
\frac{\pi}{2}-\arctan\left(\sinh(2m+1)\alpha_j\right)=(-1)^m \sum_{|k|\le m} \arctan\frac{\cos\frac{2 \pi  k}{2 m+1}}{\alpha\cos \frac{\pi  j}{2 n+1}},
$$
\end{Lemma}

\noindent{\it{Proof.}} Using properties of complex numbers we write
\begin{align*}
    \frac{\pi}{2}-\arctan\left(\sinh(2m+1)\alpha_j\right)&=\text{arg}(i)+\text{arg}\left(1-i \sinh(2m+1)\alpha_j\right)\\
    &=\text{arg}\left(\sinh(2m+1)\alpha_j+i\right)\\
    &=(-1)^m\text{arg}\left(\sinh(2m+1)\alpha_j+\sinh\tfrac{\pi i(2m+1)}{2}\right).
\end{align*}
This expression can be factorised according to the formula
$$
\sinh(2m+1)a+\sinh(2m+1)b=2^{2m}\prod_{|k|\le m}\left(\sinh a+\sinh\left(b+\frac{2\pi ik}{2m+1}\right)\right).
$$
Its validity is easy to check by standard methods: both sides are polynomials in $\sinh a$ of order $2m+1$ with leading coefficient $2^{2m}$ and zeroes $-\sinh\left(b+\frac{2\pi ik}{2m+1}\right)$, $|k|\le m $.

Thus
\begin{align*}
    \frac{\pi}{2}-\arctan\left(\sinh(2m+1)\alpha_j\right)&=(-1)^m\text{arg}\left(2^{2m}\prod_{|k|\le m}\left(\sinh \alpha_j+\sinh\left(\frac{\pi i}{2}+\frac{2\pi ik}{2m+1}\right)\right)\right)\\
    &=(-1)^m\text{arg}\left(\prod_{|k|\le m}\left(\alpha\cos\frac{\pi j}{2n+1}+i\cos\frac{2\pi k}{2m+1}\right)\right)\\
    &=(-1)^m\sum_{|k|\le m}\arctan\frac{\cos\frac{2 \pi  k}{2 m+1}}{\alpha\cos \frac{\pi  j}{2 n+1}},
\end{align*}
as required.\qed

\begin{Lemma} For $n,m\in\mathbb{N}_0$, $j\in\mathbb{Z}$, one has
$$
\sum_{|k|\le m} \arctan\frac{\cos\frac{2 \pi  k}{2 m+1}}{\alpha\cos \frac{\pi  j}{2 n+1}}=\sum_{|k|\le m} (-1)^k\arctan\frac{\cos\frac{\pi  k}{2 m+1}}{\alpha\cos \frac{\pi  j}{2 n+1}}.
$$
\end{Lemma}

\noindent{\it{Proof.}} Let $f$ be an odd function. Then
\begin{align*}
    \sum_{|k|\le m} (-1)^kf\left(\cos\frac{\pi  k}{2 m+1}\right) &= \sum_{|k|\le m} f\left(\cos\left(\frac{\pi  k}{2 m+1}-\pi k\right)\right) \\
    &=   \sum_{|k|\le m}f\left(\cos\frac{2\pi km}{2 m+1}\right)\\
    &=   \sum_{|k|\le m}f\left(\cos\frac{2\pi k}{2 m+1}\right).
\end{align*}
The last equality is explained as follows. First, note that $\cos$ has period $2\pi$. The sum $\sum_{|k|\le m}$ is over residue class mod $2m+1$. When $m>0$, the numbers $m$ and $2m+1$ are coprime. Hence, when $k$ runs over residue class mod $2m+1$, the set of numbers $km$ runs over residue class mod $2m+1$. 

To complete the proof of the lemma set ${\displaystyle{f(x)=\arctan\frac{x}{\alpha\cos \frac{\pi  j}{2 n+1}}}}$.\qed

\begin{Lemma}\label{lemma4}
$$
\sum _{|j|\le n}(-1)^j=(-1)^{n},\qquad n\in\mathbb{N}_0.
$$
\end{Lemma}

\noindent{\it{Proof.}} The sum is trivial when $n=0$. Let's assume that $n>0$. Then
\[
\pushQED{\qed} 
\sum _{|j|\le n}(-1)^j=(-1)^n\frac{1-(-1)^{2n+1}}{1-(-1)}=(-1)^n.\qedhere
\popQED
\]

Now, we are in a position to prove Theorem \ref{th1}. According to lemmas \ref{lemma1}-\ref{lemma4} we have that the LHS of equation (\ref{reciprocal_sum2}) equals
\begin{align*}
  \sum_{|j|\le n}(-1)^{n+j}\frac{1}{2}\sum_{|k|\le m} (-1)^{m+k}&\arctan\frac{\cos\frac{\pi  k}{2 m+1}}{\alpha\cos \frac{\pi  j}{2 n+1}}+  \sum_{|k|\le m} (-1)^{m+k}\frac{1}{2}\sum_{|j|\le n}(-1)^{n+j}\arctan\frac{\cos\frac{\pi  j}{2 n+1}}{\beta\cos \frac{\pi  k}{2 m+1}}\\
  &=\frac{1}{2}(-1)^{n+m}\sum_{|j|\le n}\sum_{|k|\le m}(-1)^{j+k}\,\frac{\pi}{2}\,\mathrm{sign}\left(\alpha \cos \frac{\pi  j}{2 n+1}\cos \frac{\pi  k}{2 m+1}\right)\\&=\frac{\pi}{4}(-1)^{n+m}\sum_{|j|\le n}(-1)^j\sum_{|k|\le m}(-1)^{k}=\frac{\pi}{4}.
\end{align*}

\section{Second proof of Theorem \ref{th1}}\label{proof12}

\begin{Lemma} We have the partial fractions expansion for arbitrary non-negative integer $m$:
\begin{equation}\label{part_frac}
\frac{2m+1}{\cosh \left((2m+1) \sinh ^{-1}z\right)\sqrt{z^2+1} }=\sum _{|k|\le m} \frac{(-1)^{m-k}\cos\frac{\pi k}{2m+1}}{z^2+ \cos^2\frac{\pi  k}{2m+1}}.
\end{equation}
\end{Lemma}

\noindent{\it{Proof.}} $\cosh \left((2m+1) \sinh ^{-1}z\right)\sqrt{z^2+1} $ is a polynomial in $z$ of order $2m+2$ with roots
$$
z_s=i\sin\frac{\pi(2s+1)}{2(2m+1)},\quad s=-m-1,..., m.
$$
Let us denote the LHS of \ref{part_frac} by $f(z)$. Residues of $f(z)$ at the points $z_s$ are
$$\mathrm{res}\, f(z)\big|_{z_s} =\frac{(-1)^{s}}{i(1+\delta_{s,m}+\delta_{s,-m-1})},\quad s=-m-1,..., m,$$
where $\delta_{s,r}=\begin{cases}1, s=r\\ 0, s\neq r\end{cases}$ is the Kronecker delta. It is easy to see this for $s=-m,..., m-1$. The points $z_m=i$ and $z_{-m-1}=-i$ are more tricky, in which case we write
$$
f(z)=\frac{2m+1}{{z^2+1}}\cdot\frac{\sqrt{z^2+1}}{\cosh \left((2m+1) \sinh ^{-1}z\right)},
$$
where the second multiplier does not have singularities at $z=\pm i$.

Now we can write the partial fractions expansion
$$
\frac{2m+1}{\cosh \left((2m+1) \sinh ^{-1}z\right)\sqrt{z^2+1} }=\frac{1}{i}\sum _{k=-m}^{m-1}\frac{(-1)^s}{z-i\sin\frac{\pi(2s+1)}{2(2m+1)}}+\frac{(-1)^m}{2i}\left(\frac{1}{z-i}-\frac{1}{z+i}\right).
$$
It is quite easy to bring this to the form stated in the lemma. \qed

\begin{Lemma} For arbitrary non-negative integers $n$ and $m$ we have the transformation formula
\begin{align}\label{transform}
\nonumber &\sum_{|j|\le n}(-1)^{n+j}\frac{2m+1}{\cosh\left((2m+1)\sinh^{-1}\left\{z\cos\tfrac{\pi  j}{2n+1}\right\}\right)}\frac{\cos\tfrac{\pi  j}{2n+1}}{\sqrt{1+z^2 \cos^2\tfrac{\pi  j}{2n+1}}}\\&=\frac{1}{{z}^{2}}\sum_{|k|\le m}(-1)^{m+k}\frac{2n+1}{\cosh\left((2n+1)\sinh^{-1}\left\{{z}^{-1}\cos\tfrac{\pi  k}{2m+1}\right\}\right)}\frac{\cos\tfrac{\pi  k}{2m+1}}{\sqrt{1+{z}^{-2} \cos^2\tfrac{\pi  k}{2m+1}}}.
\end{align}
\end{Lemma}

\noindent{\it{Proof.}} Multiply \ref{part_frac} by $z$ and replace $z$ with $z\cos\tfrac{\pi  j}{2n+1}$.  Then summing over $j$ one gets
\begin{align*}
    \sum_{|j|\le n}(-1)^{n+j}\frac{2m+1}{\cosh\left((2m+1)\sinh^{-1}\left\{z\cos\tfrac{\pi  j}{2n+1}\right\}\right)}\frac{z\cos\tfrac{\pi  j}{2n+1}}{\sqrt{1+z^2 \cos^2\tfrac{\pi  j}{2n+1}}}\\=\sum _{|j|\le n}\sum _{|k|\le m} (-1)^{n+m-j-k}\frac{\cos\tfrac{\pi  j}{2n+1}\cos\frac{\pi k}{2m+1}}{z\cos^2\tfrac{\pi  j}{2n+1}+{z}^{-1}\cos^2\frac{\pi  k}{2m+1}}.
\end{align*}
The RHS of this expression is symmetric under the transformation $z\to{z}^{-1}$, $n\leftrightarrow m$. As a result the LHS is also symmetric under this transformation, which implies \ref{transform}.
\qed

Integrating both sides of \ref{transform} wrt $\alpha$ using the elementary formulas $\int\frac{dx}{\sqrt{x^2+1}}=\sinh^{-1} x=\ln(\sqrt{1+x^2}+x)$, $\int\frac{dy}{\cosh y}=-2\arctan e^{-y}$ one obtains
\begin{align*}
\sum_{|j|\le n}(-1)^{n+j}\arctan& {\biggl(\!\sqrt{1+\alpha^2\cos^2\!\tfrac{\pi  j}{2n+1}}-\alpha\cos\tfrac{\pi  j}{2n+1}\!\biggr)^{2 m+1}}\\&=\sum_{|k|\le m}(-1)^{m+k}\left[\frac{\pi}{4}-\arctan\!{\biggl(\!\sqrt{1+\beta^2\cos ^2\!\tfrac{\pi  k}{2m+1}}-\beta\cos\tfrac{\pi k}{2m+1}\!\biggr)^{2 n+1}}\right].
\end{align*}
To complete the proof note that $\sum_{|k|\le m}(-1)^{m+k}$=1.

\section{Proof of Theorem \ref{th2}}\label{proof2}

Despite the fact that Theorem \ref{th2} is a particular case of Theorem \ref{th3} it is instructive to give an independent proof. Again, as we did in the previous sections, it is convenient to break the proof into several parts.
\begin{Lemma} We have the partial fractions expansion for arbitrary positive integer $m$:
\begin{equation}\label{partial_fractions}
    \frac{\sinh \left(m \tanh ^{-1}z\right)}{\sinh \left(3 m \tanh ^{-1}z\right)}\frac{1}{1-z^2}=\frac{1}{m \sqrt{3}}\sum _{k=1}^{{3 m}/{2}} \left(\frac{k}{3}\right) \frac{\tan \frac{\pi  k}{3 m}}{z^2+\tan ^2\frac{\pi  k}{3 m}}.
\end{equation}
\end{Lemma}

\noindent{\it{Proof.}} Since $\frac{\sinh t}{\sinh 3 t}=\frac{1}{2\cosh 2t+1}$, $\tanh ^{-1}z=\frac{1}{2}\ln\frac{1+z}{1-z}$, and
$$
2\cosh \left(2m \tanh ^{-1}z\right)=\left(\frac{1+z}{1-z}\right)^m+\left(\frac{1-z}{1+z}\right)^m,
$$
the LHS of (\ref{partial_fractions}) is a rational function of $z$ of the form $f(z)=\frac{(1-z^2)^{m-1}}{P_{2m}(z)}$, where $P_{2m}(z)$ is polynomial of degree exactly $2m$. This rational function has poles at $z_k=i\tan\frac{\pi k}{3m}$, $k=3r-1$ or $3r-2$ with $r=1,2,3,...,m$. Rezidues of $f(z)$ at $z_k$ are
$$
\frac{(-1)^k}{3im}\sin\frac{\pi k}{3}=\frac{1}{2im \sqrt{3}}\left(\frac{k}{3}\right).
$$
Hence, taking into account that $\left(\frac{k}{3}\right)=0$ when $k\equiv 0~(\!\!\!\!\mod{3})$
$$
\frac{\sinh \left(m \tanh ^{-1}z\right)}{\sinh \left(3 m \tanh ^{-1}z\right)}\frac{1}{1-z^2}=\frac{1}{2im \sqrt{3}}\sum _{k=1}^{3 m} \left(\frac{k}{3}\right) \frac{1}{z-i\tan\frac{\pi  k}{3 m}}.
$$
Due to $\left(\frac{3m-k}{3}\right)=-\left(\frac{k}{3}\right)$ and $\tan\frac{\pi(3m-k)}{3 m}=-\tan\frac{\pi  k}{3 m}$ this is equivalent to (\ref{partial_fractions}).\qed

\begin{Lemma} For arbitrary positive integers $n$ and $m$ we have the transformation formula
\begin{align}\label{transformation}
\nonumber m\sum _{j=1}^{{3 n}/{2}}\left(\frac{j}{3}\right)&\frac{\sinh \left(m\tanh ^{-1}\left(z\tan \frac{\pi  j}{3 n}\right)\right)}{\sinh \left(3 m\tanh ^{-1}\left(z\tan \frac{\pi  j}{3 n}\right)\right)} \frac{\tan \frac{\pi  j}{3 n}}{1-z^2\tan^2 \frac{\pi  j}{3 n}}\\&-\frac{n}{z^2}\sum _{k=1}^{{3 m}/{2}}\left(\frac{k}{3}\right)\frac{\sinh \left(n\tanh^{-1}\left(z^{-1}\tan \frac{\pi  k}{3 m}\right)\right)}{\sinh \left(3 n\tanh ^{-1}\left(z^{-1}\tan \frac{\pi  k}{3 m}\right)\right)} \frac{\tan \frac{\pi  k}{3 m}}{1-z^{-2}\tan^2 \frac{\pi  k}{3 m}}=0.
\end{align}
\end{Lemma}

\noindent{\it{Proof.}} In the previous lemma, replace $z$ with $z\tan \frac{\pi  j}{3 n}$, then multiply the resulting identity with $$\frac{z}{n}\left(\frac{j}{3}\right)\tan \frac{\pi  j}{3 n},$$ and sum wrt $j$ from $1$ to $3n/2$. It is easy to see the symmetry of the resulting double sum under the transformation $n \to m$, $m\to n$, $z \to 1/z$, from which the identity in the lemma follows. \qed

\begin{Lemma}\label{integral}
\begin{equation*}
    \sqrt{3}\int_s^\infty\frac{\sinh t}{\sinh 3 t}\, dt=\frac{\pi }{6}-\arctan\frac{\tanh s}{\sqrt{3}}=\arctan\frac{\sqrt{3}}{1+2 e^{2 s}}.
\end{equation*}
\end{Lemma}

\noindent{\it{Proof.}} The proof of this lemma is given in \cite{cais}, but we reproduce it here for the sake of completeness. Since
$$
\frac{\sinh t}{\sinh 3 t}=\frac{1}{2\cosh 2t+1}=\frac{1}{\cosh^2t}\frac{1}{3+\tanh^2 t},
$$
the integral can be easily calculated. The second equality follows from the elementary formula
$$
\arctan x-\arctan y=\arctan\frac{x-y}{1+xy}
$$
with $x=\frac{1}{\sqrt{3}}$, $y=\frac{\tanh s}{\sqrt{3}}$ and the identity
\[
\pushQED{\qed} 
\frac{\frac{1}{\sqrt{3}}-\frac{\tanh s}{\sqrt{3}}}{1+\frac{\tanh s}{3}}=\frac{\sqrt{3}}{1+2 e^s}.\qedhere
\popQED
\]

\begin{Lemma}\label{legendre}
For an odd positive integer $n$$:$
$$
\sum _{j=1}^{{3 n}/{2}} \left(\frac{j}{3}\right)=1.
$$
\end{Lemma}

\noindent{\it{Proof.}} This is obvious for $n=1$. For arbitrary odd $n$ its validity follows from the fact that the sum of Legendre symbols mod $3$ for three consecutive integers is $0$.
\qed

The formula in Theorem \ref{th2} now follows easily from these lemmas. We integrate equation (\ref{transformation}) wrt $z$ from $1/\alpha$ to $\infty$ using lemma \ref{integral}. Then assuming that $n$ and $m$ are odd we complete the proof using lemma \ref{legendre}.

\section{Proof of Theorem \ref{th3}}\label{proof3}
With the help of the elementary formula $\arctan z=\frac{1}{2}\mathrm{arg}\frac{1+iz}{1-iz}$ one can recast the first sum in \ref{parameter} in the following form
\begin{equation}\label{sum}
    \frac{1}{2} \sum _{j=0}^{\frac{n-1}{2}} \arg\frac{\left(\frac{\alpha+\tan\frac{\varphi+\pi  j}{n}}{\alpha-\tan\frac{\varphi+\pi  j}{n}}\right)^m-e^{-2i\theta}}{\left(\frac{\alpha+\tan\frac{\varphi+\pi  j}{n}}{\alpha-\tan\frac{\varphi+\pi  j}{n}}\right)^m-e^{2i\theta}}-\frac{1}{2} \sum _{j=1}^{\frac{n-1}{2}} \arg \frac{\left(\frac{\alpha-\tan\frac{\varphi-\pi  j}{n}}{\alpha+\tan\frac{\varphi-\pi  j}{n}}\right)^m-e^{-2i\theta}}{\left(\frac{\alpha-\tan\frac{\varphi-\pi  j}{n}}{\alpha+\tan\frac{\varphi-\pi  j}{n}}\right)^m-e^{2i\theta}}.
\end{equation}
In the second sum of \ref{sum}, we make the change of the index of summation $j\to n-j$, then rewrite both sums as double sums using the fact that
$x^m-1=\prod_{k=1}^{m}\big(x-e^{\frac{2\pi i k}{m}}\big)$:
\begin{equation}\label{doublesum1}
    \frac{1}{2} \sum _{j=0}^{\frac{n-1}{2}}\sum_{k=1}^m \arg\frac{\frac{\alpha+\tan\frac{\varphi+\pi  j}{n}}{\alpha-\tan\frac{\varphi+\pi  j}{n}}-e^{2i\frac{\pi k-\theta}{m}}}{\frac{\alpha+\tan\frac{\varphi+\pi  j}{n}}{\alpha-\tan\frac{\varphi+\pi  j}{n}}-e^{2i\frac{\pi k+\theta}{m}}}-\frac{1}{2} \sum_{j=\frac{n+1}{2}}^{n-1} \sum_{k=1}^m\arg
 \frac{\frac{\alpha-\tan\frac{\varphi+\pi  j}{n}}{\alpha+\tan\frac{\varphi+\pi  j}{n}}-e^{2i\frac{\pi k-\theta}{m}}}{\frac{\alpha-\tan\frac{\varphi+\pi  j}{n}}{\alpha+\tan\frac{\varphi+\pi  j}{n}}-e^{2i\frac{\pi k+\theta}{m}}}.
\end{equation}
After simple algebraic manipulation of the summands, \ref{doublesum1} becomes
\begin{align}\label{doublesum2}
    \nonumber&\frac{1}{2} \sum _{j=0}^{\frac{n-1}{2}}\sum_{k=1}^m \arg\left\{e^{-2i\frac{\theta}{m}}\frac{-i\alpha\sin\frac{\pi k-\theta}{m}+\tan\frac{\varphi+\pi  j}{n}\cos\frac{\pi k-\theta}{m}}{-i\alpha\sin\frac{\pi k+\theta}{m}+\tan\frac{\varphi+\pi  j}{n}\cos\frac{\pi k+\theta}{m}}\right\}\\&-\frac{1}{2} \sum_{j=\frac{n+1}{2}}^{n-1} \sum_{k=1}^m\arg
\left\{e^{-2i\frac{\theta}{m}}\frac{i\alpha\sin\frac{\pi k-\theta}{m}+\tan\frac{\varphi+\pi  j}{n}\cos\frac{\pi k-\theta}{m}}{i\alpha\sin\frac{\pi k+\theta}{m}+\tan\frac{\varphi+\pi  j}{n}\cos\frac{\pi k+\theta}{m}}\right\}.
\end{align}
The first sum in \ref{doublesum2} can be simplified as
$$
-\frac{n+1}{2}\theta+\frac{1}{2} \sum _{j=0}^{\frac{n-1}{2}}\sum_{k=1}^m\left(\arctan\left\{\alpha\frac{\tan\frac{\pi k+\theta}{m}}{\tan\frac{\varphi+\pi  j}{n}}\right\}-\arctan\left\{\alpha\frac{\tan\frac{\pi k-\theta}{m}}{\tan\frac{\varphi+\pi  j}{n}}\right\}\right),
$$
while the second as
$$
-\frac{n-1}{2}\theta-\frac{1}{2} \sum_{j=\frac{n+1}{2}}^{n-1} \sum_{k=1}^m\left(\arctan\left\{\alpha\frac{\tan\frac{\pi k+\theta}{m}}{\tan\frac{\varphi+\pi  j}{n}}\right\}-\arctan\left\{\alpha\frac{\tan\frac{\pi k-\theta}{m}}{\tan\frac{\varphi+\pi  j}{n}}\right\}\right).
$$
Thus, \ref{doublesum2} equals
\begin{equation}\label{doublesum3}
    -\theta+\frac{1}{2} \sum _{j=0}^{n-1}\sum_{k=1}^m\arctan\left\{\alpha\frac{\tan\frac{\pi k+\theta}{m}}{\tan\frac{\varphi+\pi  j}{n}}\right\}-\frac{1}{2} \sum _{j=0}^{n-1}\sum_{k=1}^m\arctan\left\{\alpha\frac{\tan\frac{\pi k-\theta}{m}}{\tan\frac{\varphi+\pi  j}{n}}\right\}.
\end{equation}
We make the change of the summation variable $k\to m-k$ in the second sum of \ref{doublesum3}
$$
-\theta+\frac{1}{2} \sum _{j=0}^{n-1}\sum_{k=1}^m\arctan\left\{\alpha\frac{\tan\frac{\pi k+\theta}{m}}{\tan\frac{\varphi+\pi  j}{n}}\right\}+\frac{1}{2} \sum _{j=0}^{n-1}\sum_{k=0}^{m-1}\arctan\left\{\alpha\frac{\tan\frac{\pi k+\theta}{m}}{\tan\frac{\varphi+\pi  j}{n}}\right\}.
$$
Because of the equivalences $\sum_{j=0}^{n-1}=\sum _{j=1}^{n}$, $\sum_{k=0}^{m-1}=\sum_{k=1}^m$ we finally get the following symmetric form for the first sum in \ref{parameter}
$$
-\theta+\sum _{j=1}^{n}\sum_{k=1}^m\arctan\left\{\alpha\frac{\tan\frac{\pi k+\theta}{m}}{\tan\frac{\varphi+\pi  j}{n}}\right\}.
$$
Similarly for the second sum in \ref{parameter}
$$
-\varphi+\sum _{j=1}^{n}\sum_{k=1}^m\arctan\left\{\beta\frac{\tan\frac{\varphi+\pi  j}{n}}{\tan\frac{\pi k+\theta}{m}}\right\}.
$$
Hence, because of the elementary formula $\arctan x+\arctan x^{-1}=\frac{\pi}{2}\,\mathrm{sgn}\, x$ the LHS of \ref{parameter} equals
$$
-\theta-\varphi+\frac{\pi}{2}\sum _{j=1}^{n}\mathrm{sgn}\,\left( \tan\frac{\varphi+\pi  j}{n}\right)\sum_{k=1}^m\mathrm{sgn}\left( \tan\frac{\pi k+\theta}{m}\right).
$$
When $\varphi\in(0,\pi/2)$, the expression $\tan\frac{\varphi+\pi  j}{n}$ takes negative values at $\frac{n-1}{2}$ points $j=\frac{n+1}{2},...,n-1$ and positive values at the rest $\frac{n+1}{2}$ points. This means
$$
\sum _{j=1}^{n}\mathrm{sgn}\,\left( \tan\frac{\varphi+\pi  j}{n}\right)=1.
$$
Thus the LHS of \ref{parameter} equals $\frac{\pi}{2}-\theta-\varphi$, as required.

\section{Other reciprocal relations}\label{discretedirichlet}
In our previous paper \cite{nicholson}, we have found many relations of the form $P(n,m)=P(m,n)$ for finite products of trigonometric functions. However, the identity in Theorem $1$ is of the type $S(n,m)+S(m,n)=C$, where $C$ is independent of $n$ and $m$. There is simple method to find other relations of this type. It is based on the solution of Dirichlet problem on a finite rectangular grid \cite{phillips}. For example
\begin{align}\label{dirichlet}
   m \sum _{j=1}^n (-1)^{j}\cot \frac{\pi  j}{2 n}\,\frac{ \sinh y \alpha_j}{ \sinh m \alpha_j}\sin\frac{\pi  j x}{n}+n\sum _{k=1}^m (-1)^{k}\cot \frac{\pi  k}{2 m}\, \frac{ \sinh x \beta_k}{\sinh n\beta_k}\sin\frac{\pi  ky}{m}
=-xy,
\end{align}
where $1\le x\le n$,$\ 1\le y\le m$ are integers and 
\begin{equation}\label{conditions}
\cos\frac{\pi j}{n}+\cosh\alpha_j=2,~\cos\frac{\pi k}{m}+\cosh\beta_k=2 \qquad (1\le j\le n,~1\le k\le m).
\end{equation}
In particular, when $x=y$, $n=m$ this gives the closed form summation
\begin{equation}\label{closed_form3}
    \sum _{j=1}^n (-1)^{j}\cot \frac{\pi  j}{2 n}\,\frac{ \sinh x \alpha_j}{ \sinh n \alpha_j}\sin\frac{\pi  j x}{n}=-\frac{x^2}{n},\qquad \sinh\frac{\alpha_j}{2}=\sin\frac{\pi j}{2n}.
\end{equation}

Laplace operator on a finite rectangular grid is defined as
$$
\Delta f(x,y,k)=f(x-1,y)+f(x+1,y)+f(x,y-1)+f(x,y+1)-4f(x,y).
$$
One can see that the RHS of (\ref{dirichlet}) satisfies the discrete Laplace equation 
\begin{equation*}\label{laplace}
    \Delta f(x,y)=0, \quad (0<\le x\le n, 0\le y\le m)
\end{equation*}
on a rectangular grid of size $n\times m$. Also $-xy=f_1(x,y)+f_2(x,y)$, where $f_1(x,y)$ and $f_2(x,y)$ are solutions of the Laplace equation with boundary conditions
\begin{equation}\label{dirichlet1}
  \begin{cases}
f_1(0,y)=f_1(n,y)=0,~ 0\le y\le m,\\
f_1(x,0)=0, f_1(x,m)=xm,~ 0\le x\le n,
\end{cases}  
\end{equation}
\begin{equation}\label{dirichlet2}
  \begin{cases}
f_2(x,0)=f_2(x,m)=0,~ 0\le x\le n,\\
f_2(0,y)=0, f_2(n,y)=ny,~ 0\le y\le m.
\end{cases}  
\end{equation}
Partial solutions of Laplace equation corresponding to boundary conditions (\ref{dirichlet1}) and (\ref{dirichlet2}) are given by, respectively
$$
u_j^{(1)}(x,y)=\sin\frac{\pi j x}{n}\sinh{y\alpha_j}, \quad(1\le j\le n).
$$
$$
u_k^{(2)}(x,y)=\sin\frac{\pi ky}{m}\sinh{x\beta_k}, \quad(1\le k\le m).
$$

In fact this method is quite well known and there are many examples in electrodynamics and heat conduction problems in physics (e.g., \cite{courant}).

One could generalize (\ref{dirichlet}) to include one continuous parameter $\alpha$ by requiring that $\alpha_j$ and $\beta_k$ be defined by
$$
\sinh\frac{\alpha_j}{2}=\alpha\sin\frac{\pi j}{2n},~ \sinh\frac{\beta_k}{2}=\frac{1}{\alpha}\sin\frac{\pi k}{2m}, \qquad (1\le j\le n, 1\le k\le m)
$$
instead of (\ref{conditions}). However to obtain a closed form summation we would need $\alpha=1$, so this does not generalize (\ref{closed_form3}).

\end{document}